\documentclass[12pt,letterpaper,twoside]{article}

\usepackage{amsmath,amsthm,amsfonts,amssymb}
\usepackage{times}
\usepackage{inputenc}
\usepackage{amsmath}
\usepackage{nicefrac}
\usepackage{amsthm}
\usepackage{amsfonts}
\usepackage{verbatim}
\usepackage{color}
\usepackage{bbm}
\usepackage{fancyhdr}
\usepackage[arrow, matrix, curve]{xy}

\newcommand{\pdftitle}{An exponential estimate for Hilbert space-valued Ornstein--Uhlenbeck processes}

\newcommand{\pdfauthor}{Lukas Wresch}

\author{\pdfauthor\\\small{Faculty of Mathematics, Bielefeld University, Germany, E-mail: \href{mailto:wresch@math.uni-bielefeld.de}{wresch@math.uni-bielefeld.de}}}

\title{\pdftitle}

\usepackage[plainpages=false,pdfpagelabels=true,bookmarks=true,pdfauthor={\pdfauthor},pdftitle={\pdftitle}]{hyperref}

\newtheoremstyle{Theoremstyle}     
                        {1.5em}    
                        {2.5em}    
                        {}         
                        {}         
                        {\bfseries}
                        {}         
                        {\newline} 
                        {\raisebox{0.6em}{\thmname{#1}\thmnumber{ #2}\thmnote{ (#3)}}}

\newcounter{satz}[section]
\renewcommand{\thesatz}{\thesection.\arabic{satz}}

\newcounter{ref}[satz]

\makeatletter
\def\HyPsd@CatcodeWarning#1{}

\newcommand{\customlabel}[1]{%
	 \stepcounter{ref}%
   \protected@write \@auxout{}{\string\newlabel{#1}{{\thesatz.\arabic{ref}}{\thepage}{\thesatz.\arabic{ref}}{#1}{}}}%
	\hypertarget{#1}{}%
}
	
\newcommand{\customtoc}[1]{%
  \addcontentsline{toc}{subsection}{\numberline{\arabic{section}} #1 \arabic{section}.\arabic{satz}}
	
\newcommand{\todo}[1]{{\color{red} TODO: #1} \normalfont}
}

\makeatother

\theoremstyle{Theoremstyle}
\newtheorem{lem}[satz]{Lemma}
\newtheorem{pro}[satz]{Proposition}
\newtheorem{cor}[satz]{Corollary}

\newtheorem{thm}[satz]{Theorem}
\newtheorem*{prof}{Proof}

\newcommand{\R}{\mathbb{R}}
\newcommand{\N}{\mathbb{N}}

\renewcommand{\P}{\mathbb{P}}

\newcommand{\eps}{\varepsilon}
\newcommand{\E}{\mathbb E}

\renewcommand{\d}{\,\mathrm{d}}

\newcommand{\bbm}{\mathbbm}
\newcommand{\ca}{\mathcal}
\newcommand{\mrm}{\mathrm}

\newcommand{\op}[1]{\operatorname{#1}}
								
\newcommand{\rev}{\overset{\leftarrow}}

\newcommand{\la}{\langle}
\newcommand{\ra}{\rangle}

\newcommand{\bfrac}[2]{\genfrac{}{}{0pt}{}{#1}{#2}}

\begin{document}

\maketitle

\begin{abstract}
\noindent
Let $Z$ be a $H$-valued Ornstein--Uhlenbeck process, $b\colon[0,1]\times H \rightarrow H$ and $h\colon[0,1] \rightarrow H$ be a bounded, Borel measurable functions with $\|b\|_\infty \leq 1$ then\\$\E \exp \alpha \left| \int\limits_0^1 b(t, Z_t + h(t)) - b(t, Z_t) \d t \right|_{H}^2  \leq  C$ holds, where the constant $C$ is an absolute constant and $\alpha>0$ depends only on the eigenvalues of the drift term of $Z$ and $\|h\|_\infty$, the norm of $h$, in an explicit way. Using this we furthermore prove a concentration of measure result and estimate the moments of the above integral.
\end{abstract}

\section{Introduction}

Let $H$ be a separable Hilbert space over $\R$ with an orthonormal basis $(e_n)_{n\in\N}$. Let\linebreak$(\Omega, \ca F, (\ca F_t)_{t\in[0,\infty[}, \P)$ be a filtered stochastic basis with sigma-algebra $\ca F$, a right-continuous, normal filtration $\ca F_t \subseteq \ca F$ and a probability measure $\P$ such that there is a cylindrical $\ca F_t$-Brownian motion $(B_t)_{t\in[0,\infty[}$ taking values in $\R^\N$ which is $\ca F / \ca B\left( \ca C([0,\infty[, \R^\N) \right)$ measurable, where $\ca B\left( \ca C([0,\infty[, \R^\N) \right)$ denotes the Borel sigma-algebra. We define the Wiener measure $\ca W := B(\P) := \P [ B^{-1} ]$ as the image measure of $B$ under $\P$. Let $A \colon D(A) \longrightarrow H$ be a positive definite, self-adjoint, closed, densely defined operator such that $A^{-1}$ is trace-class and
\[ \hspace{39mm}  A e_n = \lambda_n e_n, \qquad  \lambda_n > 0, \ \forall n \in\N . \]
This implies that

\customlabel{SUPLAMBDA}
\[ \sum\limits_{n\in\N} \lambda_n^{-1} =: \Lambda < \infty . \tag{\ref{SUPLAMBDA}} \]

By fixing the basis $(e_n)_{n\in\N}$ we identify $H$ with $\ell^2$, so that $H \cong \ell^2 \subseteq \R^\N$. Let $(Z_t^A)_{t\in[0,\infty[}$ be an $H$-valued (actually $\ell^2$-valued) Ornstein--Uhlenbeck process which has continuous sample paths with so-called drift term $A$, i.e.~a strong solution to
\[ \mathrm d Z_t^A = - A Z_t^A \mathrm dt + \mathrm dB_t \]
with initial condition $Z_0^A=0$.
Furthermore, we define the Ornstein--Uhlenbeck measure $\P_A$ as
\[ \P_A[F] := \P\left[ \left(Z^A\right)^{-1}(F) \right], \qquad \forall F\in \ca B( \ca C([0,\infty[, \ell^2) ). \]

In this article we show that there exist a absolute constant $C$ and an $\alpha_A>0$ such that

\customlabel{INEQ1}
\[ \E \exp \left( \frac{\alpha_A}{\|h\|_\infty^2} \left| \int\limits_0^1 b(t, Z_t^A + h(t)) - b(t, Z_t^A) \d t \right|_{H}^2 \right)  \leq  C   \tag{\ref{INEQ1}}  \]

holds, where $\alpha_A$ only depends on the eigenvalues of the drift term $A$ of $Z^A$. By approximating $b$ with smooth functions and using the Fundamental Theorem of Calculus it suffices to prove that

\customlabel{INEQ2}
\[ \E \exp \left( \alpha_A \left| \int\limits_0^1 b'(t, Z_t^A) \d t \right|_{H}^2 \right) \leq C   \tag{\ref{INEQ2}} \]

holds for one-dimensional Ornstein--Uhlenbeck processes. In order to prove \eqref{INEQ2} we follow \cite{Sha14} and split
\[ \int\limits_0^1 b'(t, Z_t^A) \d t = \int\limits_0^1 b(t, Z_t^A) \d^* t  - \int\limits_0^1 b(t, Z_t^A) \d t \]
into a forward and backward integral. Using explicit knowledge of the time-reversed process $\rev{Z^A}$ the estimate \eqref{INEQ2} is proved in Proposition \ref{THM1} and extended in Lemma \ref{THM1C} to infinite dimensions. The main result is then proved in Theorem \ref{THM1B}. Furthermore, we show in Corollaries \ref{THM2} and \ref{THM4} two applications of this result.

\section{Exponential estimate}

\begin{pro}
\label{THM1}
\customtoc{Proposition}

There exists and absolute constant $C\in\R$ and a non-increasing map
\[ \alpha \colon ]0,\infty[ \longrightarrow ]0,\infty[ \]
\[ \lambda \longmapsto \alpha_\lambda \]
with
\[ \alpha_\lambda e^{2\lambda} \lambda^{-1} \geq \frac{e}{1152}, \qquad \forall \lambda>0 . \]
such that for all one-dimensional Ornstein--Uhlenbeck processes $(Z_t^\lambda)_{t\in [0,\infty[}$ with drift term $\lambda > 0$, i.e.
\[ \begin{cases} \mrm d Z_t^\lambda = - \lambda Z_t^\lambda \mrm dt  +  \mrm dB_t,\\ \ \ \hspace{0.3mm} Z_0^\lambda = 0 . \end{cases} \]
and for all Borel measurable functions $b \colon [0,1]\times \R \longrightarrow H$, which are in the second component twice continuously differentiable with
\[ \| b \|_\infty := \sup\limits_{t\in[0,1],x\in \R} |b(t,x)|_{H} \leq 1. \]
The following inequality
\[ \E \exp \left( \alpha_\lambda \left| \int\limits_0^1 b'(t, Z_t^\lambda) \d t \right|_{H}^2 \right) \leq C \leq 3 \]
holds, where $b'$ denotes the first derivative of $b$ w.r.t.~the second component $x$.
\end{pro}

\begin{prof}

Let $(Z_t^\lambda)_{t\in [0,\infty[}$ be a one-dimensional Ornstein--Uhlenbeck process, i.e.~a strong solution to
\[ \mrm d Z_t^\lambda = - \lambda Z_t^\lambda \mrm dt  +  \mrm dB_t, \]
where $\lambda > 0$, $Z_0^\lambda=0$ and let $b \colon [0,1]\times \R \longrightarrow H$ be as in the assertion. Define
\[ \hspace{16mm}  Y_s := b(s, Z_s^\lambda), \qquad \forall s\in[0,1] \]
and denote by $(Y^n)_{n\in\N}$ the components of $Y$. Then by \cite[Remark 2.5]{BJ97} we have for every $n\in\N$
\[ \la Y^n,Z^\lambda \ra_1 = \int\limits_0^1 b_n'(s, Z_s^\lambda) \d \la Z^\lambda \ra_s = \int\limits_0^1 b_n'(s, Z_s^\lambda) \d s, \]
where $b_n$ is the $n$-th component of $b$ and the quadratic covariation $\la Y^n,Z^\lambda \ra_t$ is the uniform in probability limit of
\[ \sum\limits_{\bfrac{t_i, t_{i+1} \in \op{D}_m}{0 \leq t_{i} \leq t}} \left[ Y^n_{t_{i+1}} - Y^n_{t_i} \right] \cdot \left[ Z_{t_{i+1}}^\lambda - Z_{t_i}^\lambda \right]. \]
Moreover, applying \cite[Corollary 2.3]{BJ97} results in
\customlabel{TH1-DECOMPOSITION}
\[ \int\limits_0^1 b_n'(s, Z_s^\lambda) \d s = \la Y^n,Z^\lambda \ra_1 = \int\limits_0^1 Y^n_s \d^* Z_s^\lambda  -  \int\limits_0^1 Y^n_s \d Z_s^\lambda,  \tag{\ref{TH1-DECOMPOSITION}}  \]
where the backward integral is defined as
\customlabel{TH1-BACKWARD-INTEGRAL}
\[ \int\limits_0^t Y^n_s \d^* Z_s^\lambda  :=  - \int\limits_{1-t}^1 \rev{Y^n_s} \d \rev{Z_s^\lambda}, \qquad \forall t\in[0,1]  \tag{\ref{TH1-BACKWARD-INTEGRAL}} \]
and
\[ \hspace{16mm} \rev X_s := X_{1-s}, \qquad \hspace{14mm} \forall s\in[0,1] \]
denotes the time-reversal of a generic stochastic process $X$. Since \eqref{TH1-DECOMPOSITION} holds for all components $n\in\N$ we also have

\customlabel{TH1EQ1}
\[ \int\limits_0^1 b'(s, Z_s^\lambda) \d s = \la Y,Z^\lambda \ra_1 = \int\limits_0^1 Y_s \d^* Z_s^\lambda  -  \int\limits_0^1 Y_s \d Z_s^\lambda,  \tag{\ref{TH1EQ1}} \]
where $\la Y,Z^\lambda \ra$ is defined as $(\la Y^n,Z^\lambda \ra)_{n\in\N}$.\\

$Z^\lambda$ is an It\^{o} diffusion process with generator
\[ L_t  = a(t,x) \nabla_x + \frac12 \sigma(t,x) \Delta_x = -\lambda x \nabla_x + \frac12 \Delta_x . \]
I.e.~$a(t,x) = -\lambda x$ and $\sigma(t,x) = 1$. The probability density of $Z_t^\lambda$ w.r.t.~Lebesgue measure is
\[ p_t(x) = \sqrt\frac\lambda{\pi (1-e^{-2\lambda t})} e^{- \lambda x^2 / (1 - e^{-2\lambda t})}. \]
Observe that $a$ and $\sigma$ fulfill the conditions of \cite[Theorem 2.3]{MNS89}, hence, the drift term $\rev a$ and diffusion term $\rev \sigma$ of the generator $\rev{L}_t$ of the time-reversed process $\rev{Z^\lambda}$, is given by
\[ \rev a(t,x) = - a(1-t,x) + \frac1{p_{1-t}(x)} \nabla_x \left( \sigma(1-t,x) p_{1-t}(x) \right) = \left( \lambda - \frac{2\lambda}{1 - e^{2\lambda (t-1)}} \right) x \]
and
\[ \rev \sigma(t,x) = \sigma(1-t, x) = 1. \]
Therefore (see \cite[Remark 2.4]{BR07}), we obtain

\customlabel{TH1EQ3}
\[ \rev{Z_t^\lambda} = \rev{Z_0^\lambda}  +  \rev W_t  +  \int\limits_0^t \rev{Z_s^\lambda} \left( \lambda - \frac{2\lambda}{1 - e^{2\lambda (s-1)}} \right) \d s,  \tag{\ref{TH1EQ3}} \]

where $\rev W_t$ is a new Brownian motion defined by this equation. Set
\[ \ca G_t^0 := \sigma\left( \rev W_s - \rev W_t, t \leq s \leq 1 \right) \]
and let $\tilde{\ca G_t}$ be the completion of $\ca G_t^0$. Define
\[ \hspace{-17mm}  \ca G_t := \sigma\left( \tilde{\ca G}_{1-t} \cup \sigma(Z_1^\lambda) \right) \]
then $\rev W_t$ is a $\ca G_t$-Brownian motion (see \cite{Par86}). In conclusion we have by combining Equation \eqref{TH1EQ1} with \eqref{TH1-BACKWARD-INTEGRAL}
\[ - \int\limits_0^1 b'(s, Z_s^\lambda) \d s  =  \int\limits_0^1 b(1-s, \rev{Z_s^\lambda}) \d \rev{Z_s^\lambda}  +  \int\limits_0^1 b(s, Z_s^\lambda) \d Z_s^\lambda. \]
By plugging in \eqref{TH1EQ3} this is equal to
\[ \underbrace{\int\limits_0^1 b(1-s, \rev{Z_s^\lambda}) \d \rev W_s}_{=:\ \!I_1}  +  \underbrace{\int\limits_0^1 b(1-s, \rev{Z_s^\lambda}) \rev{Z_s^\lambda} \left( \lambda - \frac{2\lambda}{1 - e^{2\lambda (s-1)}} \right) \d s}_{=:\ \!I_2}  +  \underbrace{\int\limits_0^1 b(s, Z_s^\lambda) \d Z_s^\lambda}_{=:\ \!I_3} \]
\[ = I_1  +  I_2  +  I_3 =: I . \]
Observe that by \eqref{TH1EQ3} and the Yamada--Watanabe Theorem (see \cite[Theorem 2.1]{RSZ08}) $\rev{Z_t^\lambda}$ is a strong solution of an SDE driven by the noise $\rev W_t$, hence, $\rev{Z_t^\lambda}$ is $\ca G_t$-measurable so that the stochastic integral $I_1$ makes sense. In conclusion we get

\customlabel{TH1EQ10}
\[ \E \exp \left( \alpha_\lambda \left| \int\limits_0^1 b'(t, Z_t^\lambda) \d t \right|_H^2 \right) = \E \exp (\alpha_\lambda |I|_H^2) = \E \exp (\alpha_\lambda |I_1+I_2+I_3|_H^2)  ,  \tag{\ref{TH1EQ10}} \]

for $\alpha_\lambda$ to be defined later. We will estimate the terms $I_1$, $I_2$ and $I_3$ separately.\\

\textbf{Estimate for $I_1$:}\\

Define
\[ \hspace{18mm}  M_t := \int\limits_0^t b(1-s, \rev{Z_s^\lambda}) \d \rev W_s, \qquad \forall t\in[0,1]. \]
Observe that $(M_t)_{t\in[0,1]}$ is a $(\ca G_t)_{t\in[0,1]}$-martingale with $M_0 = 0$.
Also note the following estimate for the quadratic variation of $M$
\[ 0 \leq \la M \ra_t \leq \int\limits_0^t \| b \|_\infty^2 \d s \leq \| b \|_\infty^2 \leq 1, \qquad \forall t\in[0,1] .  \]
In the next step we use the Burkholder--Davis--Gundy Inequality for time-continuous martingales with the optimal constant. In the celebrated paper \cite[Section 3]{Dav76} it is shown that the optimal constant in our case is the largest positive root of the Hermite polynomial of order $2k$. We refer to the appendix of \cite{Ose12} for a discussion of the asymptotic of the largest positive root. See also \cite[Appendix B]{Kho14}, where a self-contained proof of the Burkholder--Davis--Gundy inequality with asymptotically optimal constant can be found for the one-dimensional case. A proof for $H$-valued martingales can be obtained by a slight modification of \cite[Theorem B.1]{Kho14} to $\R^d$-valued martingales and by projecting $H$ onto $\R^d$. The optimal constants in different cases is discussed in the introduction of \cite{Wan91}. We have
\[ \E |I_1|_{H}^{2k} = \E |M_1|_{H}^{2k} \leq 2^{2k} (2k)^{k} \underbrace{\E |\la M \ra_1|_H^{k}}_{\leq1} \leq 2^{3k} \!\! \underbrace{k^k}_{\leq 2^{2k} k!} \! \leq 2^{5k} k!   . \]
Choosing $\alpha_1 = \frac1{64}$ we obtain
\[ \E \exp \left( \alpha_1 |I_1|_{H}^2 \right) = \E \sum\limits_{k=0}^\infty \frac{\alpha_1^k I_1^{2k}}{k!} = \sum\limits_{k=0}^\infty \frac{\alpha_1^k \E |I_1|_H^{2k}}{k!} \leq \sum\limits_{k=0}^\infty 2^{-k} = 2 =: C_1. \]

\textbf{Estimate for $I_2$:}\\

We have for any $\alpha_2^{(\lambda)} > 0$ to be specified later
\[ \hspace{-31mm}  \E \exp \alpha_2^{(\lambda)} \left| I_2 \right|_{H}^2  =  \E \exp \alpha_2^{(\lambda)} \left| \int\limits_0^1 b(1-t, \rev{Z_t^\lambda})  \rev{Z_t^\lambda} \lambda \left( 1 - \frac{2}{1 - e^{2\lambda (t-1)}} \right)  \d t \right|_{H}^2 \]
\[ \hspace{-8mm}  \leq \E \exp \alpha_2^{(\lambda)} \left| \int\limits_0^1 \right. \vphantom{\int\limits_0^1} \underbrace{|b(1-t, \rev{Z_t^\lambda})|_{H}}_{\leq1} |\rev{Z_t^\lambda}| \lambda \frac{1 + e^{2\lambda(t-1)}}{1-e^{2\lambda(t-1)}}  \d t \left. \vphantom{\int\limits_0^1} \right|^2 \]
\[ \hspace{23mm} = \E \exp \alpha_2^{(\lambda)} \left| \int\limits_0^1 \frac{|\rev{Z_t^\lambda}|}{\sqrt{e^{2\lambda(1-t)}-1}} \lambda \right. \underbrace{\left(e^{2\lambda(1-t)}-1\right) \frac{1 + e^{2\lambda(t-1)}}{1-e^{2\lambda(t-1)}}}_{=e^{-2\lambda(t-1)}+1}   \left. \vphantom{\int\limits_0^1}  \frac{\mrm d t}{\sqrt{e^{2\lambda(1-t)}-1}} \right|^2 \]
\[ \hspace{0.5mm} \leq \E \exp \alpha_2^{(\lambda)} \left| \int\limits_0^1 \frac{|\rev{Z_t^\lambda}|}{\sqrt{e^{2\lambda(1-t)}-1}} \lambda  (e^{2\lambda(1-t)}+1)   \frac{\mrm d t}{\sqrt{e^{2\lambda(1-t)}-1}} \right|^2. \]
Setting
\[ \hspace{-24mm}  D_\lambda := \int\limits_0^1 \frac{\mrm d t}{\sqrt{e^{2\lambda(1-t)}-1}} =   \frac{\arctan\left(\sqrt{e^{2\lambda}-1}\right)}{\lambda}  < \infty , \]
the above term can be written as
\[ \E \exp \alpha_2^{(\lambda)} \left| \int\limits_0^1 \frac{|\rev{Z_t^\lambda}|}{\sqrt{e^{2\lambda(1-t)}-1}} \lambda \left(e^{2\lambda(1-t)}+1\right)   D_\lambda \frac{\mrm d t}{D_\lambda \sqrt{e^{2\lambda(1-t)}-1}} \right|^2 . \]
Applying Jensen's Inequality w.r.t.~the probability measure $\frac{\mrm d t}{D_\lambda \sqrt{e^{2\lambda(1-t)}-1}}$ and the convex function $\exp \alpha_2^{(\lambda)} \left| \ \! \cdot \ \!  \right|^2$ results in the above being bounded by the following
\[ \E \int\limits_0^1 \exp \left[ \alpha_2^{(\lambda)} \left| \frac{|\rev{Z_t^\lambda}|}{\sqrt{e^{2\lambda(1-t)}-1}} \lambda \left(e^{2\lambda(1-t)}+1\right) D_\lambda \right|^2 \right] \frac{\mrm d t}{D_\lambda \sqrt{e^{2\lambda(1-t)}-1}} \]
\[ \hspace{2mm}  = \E \int\limits_0^1 \exp \left[ \alpha_2^{(\lambda)} \frac{|Z_{1-t}^\lambda|^2}{e^{2\lambda(1-t)}-1} \lambda^2 \left(e^{2\lambda(1-t)}+1\right)^2 D_\lambda^2 \right] \frac{\mrm d t}{D_\lambda \sqrt{e^{2\lambda(1-t)}-1}} . \]
Setting $\alpha_2^{(\lambda)} := \frac1{4 \lambda (e^{2\lambda}+1) D_\lambda^2}$ and applying Fubini's Theorem the above term can be estimated by

\customlabel{TH1EQ20}
\[ \int\limits_0^1 \E \exp \left( \frac14 \frac{\lambda (e^{2\lambda(1-t)}+1) |Z_{1-t}^\lambda|^2}{e^{2\lambda(1-t)}-1} \right) \frac{\mrm d t}{D_\lambda \sqrt{e^{2\lambda(1-t)}-1}}. \tag{\ref{TH1EQ20}} \]

Using \cite[Theorem 8.5.7]{Oks10} (see also Step 2 of the proof of Theorem \ref{THM1B}) we have
\[ Z_{1-t}^\lambda = (2\lambda)^{-1/2} e^{-\lambda (1-t)} \bar B_{ e^{2\lambda (1-t)} - 1 }, \]
where $\bar B$ is another Brownian motion. Plugging this into \eqref{TH1EQ20} we get the following bound for \eqref{TH1EQ20}
\[ \int\limits_0^1 \E \exp \left( \vphantom{\frac{\bar B_{ e^{2\lambda (1-t)} - 1 }^2 }{4 (e^{2\lambda(1-t)}-1)}}  \right. \frac18 \frac{ \overbrace{(e^{2\lambda(1-t)}+1)e^{-2\lambda (1-t)}}^{\leq2} \bar B_{ e^{2\lambda (1-t)} - 1}^2 }{e^{2\lambda(1-t)}-1} \left.   \vphantom{\frac{\bar B_{ e^{2\lambda (1-t)} - 1 }^2 }{4 (e^{2\lambda(1-t)}-1)}}  \right) \frac{\mrm d t}{D_\lambda \sqrt{e^{2\lambda(1-t)}-1}} \]
\[ \hspace{9mm}  \leq \int\limits_0^1 \underbrace{\E \exp \left( \vphantom{\frac{ \bar B_{e^{2\lambda (1-t)} - 1}^2  }{4 (e^{2\lambda(1-t)}-1)}} \right. \frac14 \frac{ \bar B_{e^{2\lambda (1-t)} - 1}^2  }{e^{2\lambda(1-t)}-1} \left. \vphantom{\frac{ \bar B_{e^{2\lambda (1-t)} - 1}^2  }{4 (e^{2\lambda(1-t)}-1)}} \right)}_{=\sqrt2} \frac{\mrm d t}{D_\lambda \sqrt{e^{2\lambda(1-t)}-1}} \]
\[ \hspace{-5mm}  = \sqrt{2} \underbrace{\int\limits_0^1 \frac{\mrm d t}{D_\lambda \sqrt{e^{2\lambda(1-t)}-1}}}_{=1} = \sqrt{2} =: C_2. \]

\textbf{Estimate for $I_3$:}\\

Recall that

\customlabel{TH1EQ11}
\[ \E |I_3|_{H}^{2k} = \E \left| \int\limits_0^1 b(s, Z_s^\lambda) \d Z_s^\lambda \right|_{H}^{2k}. \tag{\ref{TH1EQ11}} \]

Plugging in
\[ Z_t^\lambda = - \lambda \int\limits_0^t Z_s^\lambda \d s  +  B_t \]
into Equation \eqref{TH1EQ11} results in
\[ \E |I_3|_{H}^{2k} \leq 2^{2k} \E \left| \int\limits_0^1 b(s, Z_s^\lambda) \lambda Z_s^\lambda \d s \right|_{H}^{2k}  \! +  2^{2k} \E \left| \int\limits_0^1 b(s, Z_s^\lambda) \d B_s \right|_{H}^{2k}. \]
For the first term on the right-hand side we use Jensen's Inequality and for the second term a similar calculation as for the estimate of $I_1$ yields that the above is smaller than
\[ 2^{2k} \E \int\limits_0^1 \underbrace{\|b\|^{2k}_\infty}_{\leq 1} \lambda^{2k} |Z_s^\lambda|^{2k} \d s  \! +  2^{2k} 2^{5k} k!. \]
Using Fubini's Theorem we estimate this by
\[ 2^{2k} \lambda^{2k} \int\limits_0^1 \E |Z_s^\lambda|^{2k} \d s  +  2^{2k} 2^{5k} k!  \leq  2^{2k} \lambda^{2k} \max\limits_{s\in[0,1]} \E |Z_s^\lambda|^{2k}  +  2^{2k} 2^{5k} k!. \]
With the help of \cite[Theorem 8.5.7]{Oks10} (see also Step 2 of the proof of Theorem \ref{THM1B}) we have
\[ Z_s^\lambda = (2\lambda)^{-1/2} e^{-\lambda s} \bar B_{e^{2\lambda s} - 1}, \]
where $\bar B$ is another Brownian motion. Estimating the $2k$-moments yields
\[ \hspace{-2mm} \E |Z_s^\lambda|^{2k} = (2\lambda)^{-k} e^{-\lambda 2k s} \E \left| \bar B_{e^{2\lambda s} - 1} \right|^{2k} \]
\[ \hspace{44mm}  = (2\lambda)^{-k} \underbrace{e^{-\lambda 2k s} \left| e^{2\lambda s} - 1\right|^{k}}_{\leq1} 2^k \pi^{-1/2} \Gamma\left(k + \frac12 \right) \]
\[ \hspace{48mm}  \leq \lambda^{-k} \pi^{-1/2} \Gamma\left(k + \frac12 \right)  \leq \lambda^{-k} k!, \qquad \forall s \in [0,1]. \]
Therefore, we obtain
\[ \hspace{-16mm}  \E |I_3|_{H}^{2k} \leq 2^{2k} \lambda^{2k} \max\limits_{s\in[0,1]} \E |Z_s^\lambda|^{2k}  +  2^{2k} 2^{5k} k! \]
\[  \hspace{22.5mm}  \leq 2^{2k} \lambda^{2k} \lambda^{-k} k!  +  2^{2k} 2^{5k} k! = 2^{2k} \lambda^{k} k!  +  2^{2k} 2^{5k} k! . \]
Choosing $\alpha_3^{(\lambda)} = 2^{-6} \min\left( \lambda^{-1}, 2^{-2} \right) $ we obtain
\[ \E \exp \left( \alpha_3^{(\lambda)} |I_3|_{H}^2 \right) = \E \sum\limits_{k=0}^\infty \frac{ \left|\alpha_3^{(\lambda)}\right|^k |I_3|_{H}^{2k}}{k!} = \sum\limits_{k=0}^\infty \frac{\left|\alpha_3^{(\lambda)}\right|^k \E |I_3|_{H}^{2k}}{k!} \leq \sum\limits_{k=0}^\infty 2\cdot 2^{-k} = 4  =: C_3. \]

\textbf{Final estimate:}\\

We are now ready to plug in all previous estimates to complete the proof.
Setting
\[ \alpha_\lambda := \frac19 \min(\alpha_1, \alpha_2^{(\lambda)}, \alpha_3^{(\lambda)}) \]
we conclude
\[ \E \exp( \alpha_\lambda |I|_{H}^2 ) = \E \exp( \alpha_\lambda |I_1 + I_2 + I_3|_{H}^2 ) \leq \E \exp( 3\alpha_\lambda |I_1|_{H}^2 + 3\alpha_\lambda |I_2|_{H}^2 + 3\alpha_\lambda |I_3|_{H}^2 ) \]
\[  \hspace{-12mm}  =  \E \exp( 3\alpha_\lambda |I_1|_{H}^2 ) \exp( 3\alpha_\lambda |I_2|_{H}^2 ) \exp( 3\alpha_\lambda |I_3|_{H}^2 )  .  \]
We apply the Young Inequality to split the three terms
\[ \E \frac{\exp( 3\alpha_\lambda |I_1|_{H}^2 )^3}3 + \E \frac{\exp( 3\alpha_\lambda |I_2|_{H}^2 )^3}3 + \E \frac{\exp( 3\alpha_\lambda |I_3|_{H}^2 )^3}3  \]
and using the estimates for $I_1$, $I_2$ and $I_3$ results in the following bound
\[ \E \frac{\exp( \alpha_1 |I_1|_{H}^2 )}3 + \E \frac{\exp( \alpha_2^{(\lambda)} |I_2|_{H}^2 )}3 + \E \frac{\exp( \alpha_3^{(\lambda)} |I_3|_{H}^2 )}3 \leq \frac13 (C_1 + C_2 + C_3) = \frac{6 + \sqrt2}3 \leq 3 . \]
We still need to show that the map $\alpha$ fulfills the claimed properties.\\

\textbf{Simplification of $\alpha_\lambda$:}\\

Recall that
\[ \alpha_\lambda = \frac19 \min( \alpha_1, \alpha_2^{(\lambda)}, \alpha_3^{(\lambda)}) = \frac19 \min \left( \frac1{256}, \frac1{4\lambda (e^{2\lambda}+1) D_\lambda^2}, \frac1{64\lambda} \right) \]
and
\[ D_\lambda = \frac{\arctan\left( \sqrt{ e^{2\lambda} - 1} \right)}\lambda . \]

First, we want to prove that $\alpha_\lambda$ is the same as
\[ \frac19 \min \left( \frac1{256}, \frac1{4\lambda (e^{2\lambda}+1) D_\lambda^2} \right) . \]
I.e.~$\alpha_3^{(\lambda)}$ is always larger than $\alpha_1$ or $\alpha_2^{(\lambda)}$. Note that for $\lambda \in \ ]0,4]$ $\alpha_3^{(\lambda)}$ is obviously larger than $\alpha_1$, hence it is enough to show that $\alpha_3^{(\lambda)} \geq \alpha_2^{(\lambda)}$ for all $\lambda>4$. We have
\[ \hspace{1.5mm}  2\lambda^2 + 2\lambda - \frac{10}{3\pi} \sqrt{16} \lambda + 2 \geq 0, \qquad \hspace{14.5mm}  \forall \lambda\in\R , \]
which implies
\[ \frac{10}{3\pi} \sqrt{16} \lambda \leq 2 + 2\lambda + 2 \lambda^2 \leq \sqrt{e^{2\lambda} + 1} , \qquad \forall \lambda > 4  . \]
Reordering and using that $\arctan$ is an increasing function leads us to
\[ \hspace{-18mm}  \sqrt{16} \lambda \leq \sqrt{e^{2\lambda} + 1} \frac{3\pi}{10} = \sqrt{e^{2\lambda} + 1} \arctan\left( \sqrt{1 + \frac2{\sqrt5}} \right) \]
\[ \hspace{17.5mm}  \leq \sqrt{e^{2\lambda} + 1} \arctan\left( \sqrt{ e^2 - 1} \right) \leq \sqrt{e^{2\lambda} + 1} \arctan\left( \sqrt{ e^{2\lambda} - 1} \right) \]
for all $\lambda>4$. Therefore we obtain
\[ 16 \lambda^2 \leq \left(e^{2\lambda} + 1\right) \arctan^2 \left( \sqrt{ e^{2\lambda} - 1} \right)  , \]
which finally implies
\[ \alpha_3^{(\lambda)} = \frac{1}{64\lambda} \geq \frac{\lambda}{4 \left(e^{2\lambda} + 1\right) \arctan^2 \left( \sqrt{ e^{2\lambda} - 1} \right)} = \alpha_2^{(\lambda)} . \]
In conclusion we proved that
\[ \alpha_\lambda = \frac19 \min \left( \frac1{256}, \frac1{4\lambda (e^{2\lambda}+1) D_\lambda^2} \right) . \]

\textbf{Asymptotic behavior of $\alpha_\lambda$:}\\

Let us now analyze $\alpha_2^{(\lambda)}$. Set
\[ f(\lambda) := \alpha_2^{(\lambda)} e^{2\lambda} \lambda^{-1} = \frac{e^{2\lambda}}{4 \lambda^2 (e^{2\lambda}+1) D_\lambda^2} = \frac{e^{2\lambda}}{4 (e^{2\lambda}+1) \arctan^2(\sqrt{e^{2\lambda} - 1})} . \]
We obviously have
\[ \frac{e^{2\lambda}}{e^{2\lambda}+1} \overset{\lambda\rightarrow\infty}\longrightarrow 1 \]
and
\[ \arctan\left(\sqrt{e^{2\lambda} - 1}\right) \overset{\lambda\rightarrow\infty}\longrightarrow \frac\pi2 . \]
Therefore,
\[ f(\lambda) = \frac{e^{2\lambda}}{4 (e^{2\lambda}+1) \arctan^2(\sqrt{e^{2\lambda} - 1})} \overset{\lambda\rightarrow\infty}\longrightarrow \frac1{\pi^2} \]
holds. We want to show that $f$ is monotonically decreasing and hence the above limit is a lower bound for $f$. To this end we calculate the first derivative of $f$
\[ f'(\lambda) = - \frac{e^{2\lambda} \left( e^{2\lambda} + 1 - 2\arctan\left( \sqrt{e^{2\lambda} - 1}\right) \sqrt{e^{2\lambda} - 1} \right) } { 4 \arctan^3\left( \sqrt{e^{2\lambda} - 1}\right) \sqrt{e^{2\lambda} - 1} \left( e^{2\lambda} - 1 \right)^2 } . \]
since the denominator is clearly positive, we have to show that
\[ e^{2\lambda} + 1 - 2\arctan\left( \sqrt{e^{2\lambda} - 1}\right) \sqrt{e^{2\lambda} - 1} > 0 , \qquad \forall \lambda > 0  .  \]
Substituting $x := \sqrt{e^{2\lambda}-1}$ leads to
\customlabel{THM-X2-ARCTAN}
\[ x^2 + 2  > 2x \arctan(x) , \qquad \forall x > 0  .   \tag{\ref{THM-X2-ARCTAN}}  \]
We prove this inequality in two steps. First note that
\[ \hspace{4mm}  x^2 - \frac{10\pi}{12} x + 2 > 0 , \hspace{-0.5mm} \qquad \forall x>0 \]
holds, so that for all $x$ with $0 < x \leq 2+\sqrt3$ we have the estimate
\[ x^2 + 2 > 2x \frac{5\pi}{12} = 2x \arctan(2+\sqrt3) \geq 2x \arctan(x) \]
and, on the other hand, for $x\geq2+\sqrt3$ we obtain
\[ x^2 + 2 \geq (2+\sqrt3)x + 2 > (2+\sqrt3)x \geq \pi x = 2x \frac\pi2 > 2x \arctan(x) . \]
In conclusion \eqref{THM-X2-ARCTAN} holds, so that $f'<0$ and therefore
\[ f(\lambda) \geq \frac1{\pi^2}, \qquad \forall \lambda > 0 . \]
All together this yields
\[ \alpha_\lambda e^{2\lambda} \lambda^{-1} = \frac19 \min\left( \frac1{256} e^{2\lambda} \lambda^{-1}, \alpha_2^{(\lambda)} e^{2\lambda} \lambda^{-1}\right) \geq \frac19 \min\left( \frac{e}{128}, \frac1{\pi^2} \right) = \frac{e}{1152} . \]

\textbf{$\alpha_\lambda$ is constant on $[0,1]$:}\\

\textbf{Claim:}
\[ \alpha_2^{(\lambda)} \geq \frac1{256}, \qquad \forall\lambda \in [0,1] . \]
Let $\lambda \in [0,1]$ and set
\[ g(\lambda) := \frac{ (e^{2\lambda}+1) (e^{2\lambda}-1)}{\lambda} . \]
$g$ has the first derivative
\[ g'(\lambda) = \frac{1 - (1-4\lambda) e^{4\lambda}}{\lambda^2} . \]
We want to show that $1 - (1-4\lambda) e^{4\lambda}$ is non-negative and thus prove that $g$ is an non-decreasing function. To this end observe that
\[ (1 - 4\lambda) e^{4\lambda} \]
is a decreasing function on $[0,\infty[$, since the derivative $- 16 \lambda e^{4\lambda}$ is clearly non-positive, so that
\[  \hspace{3mm}  (1 - 4\lambda) e^{4\lambda} \leq 1 \]
holds for all $\lambda\geq0$. This leads to
\[  \hspace{26mm}  1 - (1-4 \lambda) e^{4\lambda} \geq 0 , \qquad  \hspace{7mm}  \forall \lambda\geq0 . \]
This proves that $g$ is non-decreasing. Using this we can easily conclude
\[ \max\limits_{\lambda\in[0,1]} g(\lambda) \leq g(1) = (e^2 + 1) (e^2 - 1) \leq 64 \]
and hence

\[ g(\lambda) = \frac{ (e^{2\lambda}+1) (e^{2\lambda}-1)}{\lambda} \leq \frac{256}{4}, \qquad \forall \lambda\in [0,1] .   \]

Taking the reciprocal on both sides yields

\customlabel{G-BOUND}
\[ \alpha_2^{(\lambda)} = \frac{\lambda}{ 4 (e^{2\lambda}+1) (e^{2\lambda}-1)} \geq \frac{1}{256}, \qquad \forall \lambda\in [0,1]  .  \tag{\ref{G-BOUND}} \]

Note that
\[ \arctan(x) \leq x , \qquad \forall x\in\R_+ . \]
This can be proved by calculating the Taylor-polynomial up to the first order and dropping the remainder term which is always negative on $\R_+$. Using this on our above estimate \eqref{G-BOUND} we obtain
\[  \hspace{-5mm}  \frac{\lambda}{ 4 (e^{2\lambda}+1) \arctan^2\left(\sqrt{e^{2\lambda}-1}\right)} \geq \frac1{256}, \qquad  \hspace{2mm}  \forall \lambda\in[0,1] . \]

This implies that $\alpha_\lambda$ is constant on the interval $[0,1]$.\\

\textbf{$\alpha_\lambda$ is non-increasing:}\\

By the previous part we can assume that $\lambda\geq1$. We have to show that $\alpha_2^{(\lambda)}$ is non-increasing on the interval $[1, \infty[$. We do this by showing that the derivative of $\alpha_2^{(\lambda)}$
\[ \left( \alpha_2^{(\lambda)} \right)' = - \frac{\overbrace{2\lambda}^{=:\ \!p_1} - \overbrace{\arctan\left( \sqrt{e^{2\lambda}-1} \right) \sqrt{e^{2\lambda}-1}}^{=:\ \!n_1} }                          {4 \arctan^3 \left( \sqrt{e^{2\lambda}-1} \right) \sqrt{e^{2\lambda}-1} (e^{2\lambda} + 1)^2 } \]

\[ -  \frac{\overbrace{2\lambda e^{2\lambda}}^{=:\ \!p_2} - \overbrace{\arctan\left( \sqrt{e^{2\lambda}-1} \right) \sqrt{e^{2\lambda}-1} e^{2\lambda}}^{=:\ \!n_2} + \overbrace{2\lambda \arctan\left( \sqrt{e^{2\lambda}-1} \right) \sqrt{e^{2\lambda}-1} e^{2\lambda}}^{=:\ \!p_3} }                    {4 \arctan^3 \left( \sqrt{e^{2\lambda}-1} \right) \sqrt{e^{2\lambda}-1} (e^{2\lambda} + 1)^2 } . \]
is non-positive. So, to simplify notation we have to show that

\customlabel{P-N}
\[ p_1 - n_1 + p_2 - n_2 + p_3 \geq 0 , \qquad \forall \lambda\geq1 \tag{\ref{P-N}} \]

holds.
Note that for $\lambda\geq1$
\[ p_3 - n_1 - n_2 \geq \arctan\left( \sqrt{e^{2\lambda}-1} \right) \sqrt{e^{2\lambda}-1} e^{2\lambda} \left( 2\lambda - 2 \right) \geq 0  , \]

so that \eqref{P-N} holds, which finishes the proof that $\alpha_2^{(\lambda)}$ is non-increasing on $[1,\infty[$. Together with the previous established result that $\alpha$ is constant on $[0,1]$ this completes the proof that $\alpha_\lambda$ is non-increasing on $\R_+$.

\qed
\end{prof}

\begin{lem}
\label{THM1C}
\customtoc{Lemma}

Let $(Z_t^A)_{t\in [0,\infty[}$ be a $H$-valued Ornstein--Uhlenbeck process with drift term $A$ as explained in the introduction, i.e.
\[ \begin{cases} \mrm d Z_t^A = - A Z_t^A \mrm dt  +  \mrm dB_t,\\ \ \ \hspace{0.4mm} Z_0^A = 0 . \end{cases} \]
Let $(\lambda_n)_{n\in\N}$ be the eigenvalues of $A$. Let $C\in\R$ and the map $\alpha$ be as in Proposition \ref{THM1}. Then for all Borel measurable functions $b \colon [0,1]\times H \longrightarrow H$, which are in the second component twice continuously differentiable with
\[ \| b \|_\infty := \sup\limits_{t\in[0,1],x\in \R} |b(t,x)|_{H} \leq 1 \]
we have
\[ \E \exp \left( \alpha_{\lambda_i} \left| \int\limits_0^1 \partial_{x_i} b(t, Z_t) \d t \right|_{H}^2 \right) \leq C \leq 3 \qquad \forall i\in\N , \]
where $\partial_{x_i} b$ denotes the derivative of $b$ w.r.t.~the $i$-th component of the second parameter $x$.

\end{lem}

\begin{prof}

Let
\[ Z_t^A = (Z_t^{A,(n)})_{n\in\N} \in \ell^2 \cong H \]
be the components of $(Z_t^A)_{t\in[0,\infty[}$ and $(\lambda_n)_{n\in\N}$ be the eigenvalues of $A$ w.r.t.~the basis $(e_n)_{n\in\N}$. Note that every component $Z^{A,(n)}$ is a one-dimensional Ornstein--Uhlenbeck process with drift term $\lambda_n>0$ driven by the one-dimensional Brownian motion $B^{(n)}$. Define $\tilde B^{(n)}$ by
\[ \tilde B^{(n)}_t := \int\limits_0^{\gamma^{(n)}(t)} \sqrt{c^{(n)}(s)} \d B^{(n)}_s, \qquad \forall t\in[0,1], \]
where
\[ \gamma^{(n)}(t) := (2\lambda_n)^{-1} \ln(t+1) \qquad \text{and} \qquad c^{(n)}(t) := (2\lambda_n) e^{2\lambda_n t} . \]
Observe that
\[ \left(\gamma^{(n)}(t) \right)' = \frac1{c^{(n)} \left(\gamma^{(n)}(t) \right)} \]
and, hence, by \cite[Theorem 8.5.7]{Oks10} $\left(\tilde B^{(n)}_t\right)_{t\in[0,\infty[}$ is a Brownian motion for every $n\in\N$. Therefore, the components of $Z^A$ can be written as time-transformed Brownian motions
\[ Z_t^{A,(n)}  =  (2 \lambda_n)^{-1/2} e^{-\lambda_n t} \tilde B^{(n)}_{e^{2\lambda_n t}-1}, \]
since
\[ \hspace{-11mm}  (2 \lambda_n)^{-1/2}  \tilde B^{(n)}_{e^{2\lambda_n t}-1} = (2 \lambda_n)^{-1/2} \int\limits_0^t (2\lambda_n)^{1/2} e^{\lambda_n s} \d B_s^{(n)} \]
\[  \hspace{15mm}  =  \int\limits_0^t e^{\lambda_n s} \d B_s^{(n)} = Z_t^{A,(n)} e^{\lambda_n t}. \]
Let us define the mapping
\[ \hspace{-41mm}  \varphi_A\colon \ca C( [0,\infty[, H) \longrightarrow \ca C( [0, \infty[, H) \]
\[ \hspace{22mm}  f = (f^{(n)})_{n\in\N} \longmapsto \left( t \longmapsto \left( (2\lambda_n)^{-1/2} e^{-\lambda_n t} f^{(n)}(e^{2\lambda_n t} - 1 ) \right)_{n\in\N} \right) .   \]
$\varphi_A$ is bijective and we have used that $\ca C( [0,\infty[, \R^\N) \cong \ca C( [0,\infty[, \R)^{\N}$ as topological spaces. By definition of the product topology $\varphi_A$ is continuous if and only if $\pi_n \circ \varphi_A$ is continuous for every $n\in\N$.
\[ \begin{xy}
  \xymatrix{
      \ca C( [0,\infty[, H) \ar[rrr]^{\!\!\!\!\!\!\!\!\!\!\!\!\!\!\!\!\!\!\varphi_A} \ar[rrrdd]_{\varphi_A^{(n)} := \pi_n\circ\varphi_A}  &  & &  \ca C( [0, \infty[, H) \ar[dd]^{\pi_n}  \\
			            \\
                   &  &  &  \ca C( [0, \infty[, \R)
  }
\end{xy}
\]

Here, $\pi_n$ denotes the projection to the $n$-th component. The above mapping $\varphi_A$ is continuous and, therefore, measurable w.r.t.~the Borel sigma-algebra.
Using this transformation, the Ornstein--Uhlenbeck measure $\P_A$, as defined in the introduction, can be written as
\[ \P_A[F] =  Z^A(\P)[F] = (\varphi_A \circ \tilde B) (\P) [F] = \varphi_A(\ca W)[F] , \qquad \forall F\in\ca B\left( \ca C\left( [0,\infty[, H \right)\right), \]
because of
\[ \hspace{4.5mm}  Z_t^A = \varphi_A \circ \tilde B_t .    \] Hence, we have

\customlabel{TH1BEQ220}
\[ \P_A = \varphi_A (\ca W) = \varphi_A \left( \bigotimes\limits_{n\in\N} \ca W^{(n)} \right) = \bigotimes\limits_{n\in\N} \varphi_A^{(n)} \left( \ca W^{(n)} \right),  \tag{\ref{TH1BEQ220}}  \]

where $\ca W^{(n)}$ is the projection of $\ca W$ to the $n$-th coordinate and the last equality follows from

\[ \int\limits_F \mathrm d \varphi_A\left( \bigotimes\limits_{n\in\N} \ca W^{(n)} \right)  =  \prod\limits_{n\in\N} \int\limits_{\pi_n\left(\varphi_A^{-1}(F)\right)} \!\!\!\!\!\! \mathrm d \ca W^{(n)}  =  \prod\limits_{n\in\N} \int\limits_{(\varphi_A^{(n)})^{-1}( \pi_n(F) )} \!\!\!\!\!\!\!\!\! \mathrm d \ca W^{(n)}  = \left(\bigotimes\limits_{n\in\N} \varphi_A^{(n)} \left( \ca W^{(n)} \right) \right) [F]  .  \]

Starting from the left-hand side of the assertion we have
\[ \E \exp \alpha_{\lambda_i} \left| \int\limits_0^1 \partial_{x_i} b(t, Z_t^A) \d t \right|_{H}^2 . \]
Using Equation \eqref{TH1BEQ220} we can write this as
\[    \int\limits_{\ca C([0,\infty[, \R)^{\N}}   \exp \alpha_{\lambda_i}  \left| \int\limits_0^1 \partial_{x_i} b\left(t,  ((\varphi_A^{(n)} \circ f_n) (t))_{n\in\N} \right) \d t \right|_{H}^2 \d \bigotimes\limits_{n\in\N} \ca W^{(n)}(f_n) , \]
where $(f_n)_{n\in\N}$ are the components of $f$. Using Fubini's Theorem we can perform the $i$-th integral first and obtain
\[ \int\limits_{\ca C([0,\infty[, \R)^{\N\setminus\{i\}}}  \int\limits_{\ca C([0,\infty[, \R)}   \!\!\!   \exp \alpha_{\lambda_i} \left| \int\limits_0^1 \partial_{x_i} b\left(t,  ((\varphi_A^{(n)} \circ f_n) (t))_{n\in\N} \right) \d t \right|_{H}^2   \!  \d \ca W^{(i)}(f_i)  \mrm d \bigotimes\limits_{\bfrac{n\in\N}{n\neq i}} \ca W^{(n)}(f_n) . \]
Since $\varphi_A^{(i)} \circ f_i$ is under $\ca W^{(i)}$ distributed as $Z^{A,(i)}$ under $\P$. By Proposition \ref{THM1} the inner integral is smaller than $C$, so that the entire expression is smaller than
\[  \int\limits_{\ca C([0,\infty[, \R)^{\N \setminus \{i\}}} C  \d \bigotimes\limits_{\bfrac{n\in\N}{n\neq i}} \ca W^{(n)}(f_n) = C , \]
where in the last step we used that $\ca W^{(n)}$ are probability measures.\\
\qed
\end{prof}

\begin{thm}
\label{THM1B}
\customtoc{Theorem}

Let $\ell\in\ ]0,1]$ and $(Z_t^{\ell A})_{t\in [0,\infty[}$ be an $H$-valued Ornstein--Uhlenbeck process with drift term $\ell A$, i.e.
\[ \begin{cases} \mrm d Z_t^{\ell A} = - \ell A Z_t^{\ell A} \mrm dt  +  \mrm dB_t,\\ \ \ \!\hspace{0.8mm} Z_0^{\ell A} = 0 . \end{cases} \]
There exists an absolute constant $C\in\R$ (independent of $A$ and $\ell$) such that for all Borel measurable functions $b \colon [0,1]\times H \longrightarrow H$ with
\[ \| b \|_\infty := \sup\limits_{t\in[0,1],x\in H} |b(t,x)|_{H} \leq 1 \]
and
\[ \| b \|_{\infty,A} := \sup\limits_{t\in[0,1],x\in H} \left( \sum\limits_{n\in\N} \lambda_n e^{2\lambda_n} b_n(t,x)^2 \right)^{1/2} \leq 1 . \]
This means
\[ \qquad \left( \lambda_n^{1/2} e^{\lambda_n} b_n(t,x) \right)_{n\in\N} \in \ell^2 \cong \in H, \qquad\quad \forall (t,x) \in [0,1]\times H , \]
where $b_n$ is the $n$-th component of $b$. The following inequality
\[ \E \exp \frac{\beta_A}{\|h\|_{\infty}^2} \left| \int\limits_0^1 b(t, Z_t^{\ell A} + h(t)) - b(t, Z_t^{\ell A}) \d t \right|_{H}^2 \leq C \leq 3, \]
where
\[ \beta_A := \frac14 \Lambda^{-2} \inf\limits_{n\in\N} \alpha_{\lambda_n} e^{2\lambda_n} \lambda_n^{-1} > 0 \]
holds uniformly for all bounded, measurable functions $h \colon [0,1]\longrightarrow H$ with
\[ \|h\|_{\infty} := \sup\limits_{t\in[0,1]} |h(t)|_H \in \ \!  ]0,\infty[  \]
and
\[ \qquad \sum\limits_{n\in\N} |h_n(t)|^2 \lambda_n^2 < \infty , \quad \qquad \forall t \in [0,1]  . \]
Recall that $\Lambda$ is defined in Equation \eqref{SUPLAMBDA} and the map $\alpha$ is from Proposition \ref{THM1}.
\end{thm}

\begin{prof}

\textbf{Step 1:} The case for twice continuously differentiable $b$.\\

Let $Z^{\ell A}$ be an $H$-valued Ornstein--Uhlenbeck process, $b \colon [0,1] \times H \longrightarrow H$ a bounded, Borel measurable function which is twice continuously differentiable in the second component with $\|b\|_\infty \leq 1$, $\|b\|_{\infty,A} \leq 1$ and $h \colon [0,1]\longrightarrow H$ a bounded, measurable function with $\|h\|_{\infty} \neq0$. Let $\alpha$ and $C$ be as in Proposition \ref{THM1}. Recall that $\Lambda$ is defined as
\[ \Lambda = \sum\limits_{n\in\N} \lambda_n^{-1} < \infty . \]
Note that by Proposition \ref{THM1} $\beta_A > 0$. By the Fundamental Theorem of Calculus we obtain
\[ \E \exp \frac{4\beta_A}{\|h\|_{\infty}^2} \left| \int\limits_0^1 b(t, Z_t^{\ell A} + h(t)) - b(t, Z_t^{\ell A}) \d t \right|_{H}^2 \]
\[ = \E \exp \frac{4\beta_A}{\|h\|_{\infty}^2} \left| \int\limits_0^1 \left. \vphantom{\int} b(t, Z_t^{\ell A} + \theta h(t)) \right|_{\theta=0}^{\theta=1} \d t \right|_{H}^2 \]
\[ = \E \exp \frac{4\beta_A}{\|h\|_{\infty}^2} \left| \int\limits_0^1 \int\limits_0^1 b'(t, Z_t^{\ell A} + \theta h(t)) h(t)  \d\theta \mrm dt \right|_{H}^2 \! , \]

where $b'$ denotes the Fr\'echet derivative of $b$ w.r.t.~$x$. Using Fubini's Theorem we can switch the order of integration, so that the above equals

\[ \E \exp 4\beta_A \left| \int\limits_0^1 \int\limits_0^1 \right.  b'(t, Z_t^{\ell A} + \theta h(t)) \frac{h(t)}{\|h\|_{\infty}}  \d t \mrm d\theta \left. \vphantom{\int\limits_0^1} \right|_{H}^2 \]
\[ = \E \exp 4\beta_A \left| \int\limits_0^1 \int\limits_0^1 \right.  \sum\limits_{i\in\N} \underbrace{b'(t, Z_t^{\ell A} + \theta h(t)) e_i}_{= \partial_{x_i} b(t, Z_t^{\ell A} + \theta h(t))} \frac{h_i(t)}{\|h\|_{\infty}}  \d t \mrm d\theta \left. \vphantom{\int\limits_0^1} \right|_{H}^2 \]
\[ = \E \exp 4\beta_A \left| \int\limits_0^1 \int\limits_0^1 \right.  \sum\limits_{i\in\N} \frac{h_i(t)}{\|h\|_{\infty}} \sum\limits_{j\in\N} \partial_{x_i} b_j(t, Z_t^{\ell A} + \theta h(t)) e_j  \d t \mrm d\theta \left. \vphantom{\int\limits_0^1} \right|_{H}^2 \]

\customlabel{TH2SWAP}
\[ = \E \exp 4\beta_A \left| \int\limits_0^1 \int\limits_0^1 \right.  \sum\limits_{i\in\N} \lambda_i^{-1/2} \partial_{x_i} \underbrace{\frac{h_i(t)}{\|h\|_{\infty}} \lambda_i^{1/2} \sum\limits_{j\in\N} b_j(t, Z_t^{\ell A} + \theta h(t)) e_j}_{= e^{-\lambda_i} \tilde b_{h,\theta,i}(t, Z^{\ell A}_t)}  \d t \mrm d\theta \left. \vphantom{\int\limits_0^1} \right|_{H}^2 ,  \tag{\ref{TH2SWAP}} \]
where
\[ \tilde b_{h,\theta,i}(t, x) := e^{\lambda_i} \frac{h_i(t)}{\|h\|_{\infty}} \lambda_i^{1/2} \sum\limits_{j\in\N} b_j(t, x + \theta h(t)) e_j . \]

Note that $\|\tilde b_{h,\theta,i}\|_\infty \leq 1$ because for all $(t,x)\in [0,1]\times H$ we have
\[ | \tilde b_{h,\theta,i}(t,x) |_H = \underbrace{\frac{|h_i(t)|}{\|h\|_{\infty}}}_{\leq1} \lambda_i^{1/2} e^{\lambda_i} \left| \sum\limits_{j\in\N} b_j(t, x + \theta h(t)) e_j \right|_H \]
\[ \leq \lambda_i^{1/2} e^{\lambda_i} \left( \sum\limits_{j\in\N} \lambda_j^{-1} e^{-2\lambda_j} \lambda_j e^{2\lambda_j} b_j(t, x + \theta h(t))^2 \right)^{1/2} \]
\[ \leq \underbrace{\lambda_i^{1/2} e^{\lambda_i} \sup\limits_{j\in\N} \lambda_j^{-1/2} e^{-\lambda_j}}_{\leq 1} \underbrace{\left( \sum\limits_{j\in\N} \lambda_j e^{2\lambda_j} b_j(t, x + \theta h(t))^2 \right)^{1/2}}_{\leq1} \leq 1 . \]

Using Jensen's Inequality and again Fubini's Theorem the Expression \eqref{TH2SWAP} is bounded from above by
\[ \int\limits_0^1 \E \exp 4\beta_A \left| \sum\limits_{i\in\N} \lambda_i^{-1/2} \int\limits_0^1  e^{-\lambda_i} \partial_{x_i} \tilde b_{h,\theta,i}(t, Z_t^{\ell A})  \d t \right|_{H}^2 \mrm d\theta . \]
Applying H\"older Inequality we can split the sum and estimate this from above by
\[  \hspace{-3mm}  \int\limits_0^1 \E \exp 4\beta_A \underbrace{\sum\limits_{i\in\N} \lambda_i^{-1}}_{=\Lambda}  \sum\limits_{i\in\N} \left| \int\limits_0^1 e^{-\lambda_i} \partial_{x_i} \tilde b_{h,\theta,i}(t, Z_t^{\ell A})  \d t \right|_{H}^2 \mrm d\theta \]
\[ = \int\limits_0^1 \E \exp 4\beta_A \Lambda \sum\limits_{i\in\N} \left| \int\limits_0^1 e^{-\lambda_i} \partial_{x_i} \tilde b_{h,\theta,i}(t, Z_t^{\ell A})  \d t \right|_{H}^2 \mrm d\theta \]
\[ \hspace{1mm} = \int\limits_0^1 \E \prod\limits_{i\in\N} \exp 4\beta_A \Lambda \left| \int\limits_0^1 e^{-\lambda_i} \partial_{x_i} \tilde b_{h,\theta,i}(t, Z_t^{\ell A})  \d t \right|_{H}^2 \mrm d\theta . \]
Young's Inequality with $p_i := \lambda_i \Lambda$ leads us to the upper bound
\[ \int\limits_0^1 \E \sum\limits_{i\in\N} \frac1{p_i} \exp 4\beta_A \Lambda p_i \left| \int\limits_0^1 e^{-\lambda_i} \partial_{x_i} \tilde b_{h,\theta,i}(t, Z_t^{\ell A})  \d t \right|_{H}^2 \mrm d\theta . \]

\customlabel{TH2FINAL}
\[ \hspace{-5mm}  = \int\limits_0^1 \sum\limits_{i\in\N} \frac1{p_i} \E \exp 4\beta_A \Lambda^2 \lambda_i \left| \int\limits_0^1 e^{-\lambda_i} \partial_{x_i} \tilde b_{h,\theta,i}(t, Z_t^{\ell A})  \d t \right|_{H}^2 \mrm d\theta  \tag{\ref{TH2FINAL}} . \]

Recall that
\[ \beta_A = \frac14 \Lambda^{-2} \inf\limits_{n\in\N} \alpha_{\lambda_n} e^{2\lambda_n} \lambda_n^{-1} , \]
hence, we can estimate \eqref{TH2FINAL} from above by
\[ \int\limits_0^1 \sum\limits_{i\in\N} \frac1{p_i} \E \exp \alpha_{\lambda_i} e^{2\lambda_i} \left| \int\limits_0^1 e^{-\lambda_i} \partial_{x_i} \tilde b_{h,\theta,i}(t, Z_t^{\ell A})  \d t \right|_{H}^2 \mrm d\theta . \]
\[ = \int\limits_0^1 \sum\limits_{i\in\N}\frac1{p_i} \E \exp \alpha_{\lambda_i} \left| \int\limits_0^1 \partial_{x_i} \tilde b_{h,\theta,i}(t, Z_t^{\ell A})  \d t \right|_{H}^2 \mrm d\theta . \]
Since $\ell\in \ ]0,1]$ and $\alpha$ is non-increasing by Proposition \ref{THM1} the above is smaller than
\[ \int\limits_0^1 \sum\limits_{i\in\N} \frac1{p_i} \E \exp \alpha_{\ell \lambda_i} \left| \int\limits_0^1 \partial_{x_i} \tilde b_{h,\theta,i}(t, Z_t^{\ell A})  \d t \right|_{H}^2 \mrm d\theta . \]
Applying Lemma \ref{THM1C} for every $\theta\in[0,1]$ and $i\in\N$ results in the estimate
\[ \int\limits_0^1 \underbrace{\sum\limits_{i\in\N} \frac1{p_i}}_{=1} C \d\theta = C . \]

\textbf{Step 2:} The general case: Non-smooth $b$.\\

Let $b \colon [0,1] \times H \longrightarrow H$ be a bounded, Borel measurable function with $\|b\|_\infty \leq 1$, $\|b\|_{\infty,A} \leq 1$ and $h \colon [0,1]\longrightarrow H$ a bounded, Borel measurable function with $0 \neq \|h\|_{\infty} < \infty$ and
\[ \sum\limits_{n\in\N} |h_n(t)|^2 \lambda_n^2 < \infty \qquad \forall t \in [0,1]  . \]
Let $\beta_A$ and $C$ be the constants from Step $1$. Set $\eps := \exp \frac{- 64 \beta_A}{\|h\|_{\infty}^2}$ as well as
\[ \mu_0 := \mrm dt \otimes Z_t^{\ell A} (\P) , \]
\[ \mu_h := \mrm dt \otimes (Z_t^{\ell A} + h(t)) (\P) . \]

Note that the measure $Z^{\ell A}(\P)$ is equivalent to the invariant measure $N(0, \frac1{2\ell} A^{-1})$ due to \cite[Theorem 11.13]{DZ92} and analogously $(Z_t^{\ell A} + h(t))(\P)$ to $N(h(t), \frac1{2\ell} A^{-1})$. Furthermore, $h(t)$ is in the domain of $A$ for every $t\in[0,1]$ because of

\[ \sum\limits_{n\in\N} \la h(t), e_n \ra^2 \lambda_n^2 \leq \sum\limits_{n\in\N} |h_n(t)|^2 \lambda_n^2 < \infty . \]

We set
\[ g(t) := 2 \ell A h(t) . \]
Observe that $g(t)\in H$ for every $t\in[0,1]$ because of

\[ |g(t)|_H^2  = 4 \ell^2 \sum\limits_{n\in\N} \lambda_n^2 |h_n(t)|^2 < \infty . \]

Hence, \cite[Corollary 2.4.3]{Bog98} is applicable i.e.~$N(0, \frac1{2\ell} A^{-1})$ and $(Z_t^{\ell A} + h(t))(\P)$ are equivalent measures. By the Radon--Nikodym Theorem there exist a density $\rho$ so that
\[ \frac{\mrm d \mu_h}{\mrm d \mu_0} = \rho . \]

Furthermore, there exists $\delta > 0$ such that
\customlabel{TH2-RHO}
\[ \int\limits_A \rho  \d \mu_0(t,x) \leq \frac\eps2 ,  \tag{\ref{TH2-RHO}}  \]
for all measurable sets $A \subseteq [0,1] \times H$ with $\mu_0[A] \leq \delta$. Set
\customlabel{TH2-BARDELTA}
\[ \bar\delta := \min\left( \delta, \frac\eps2 \right).  \tag{\ref{TH2-BARDELTA}}  \]

By Lusin's Theorem (see \cite[Theorem 1.3.28]{Tao11}) there exist a closed set $K \subseteq [0,1] \times H$ with $\mu_0[K] \geq 1 - \bar\delta$ such that the restriction
\[ b\mid_K\colon K \longrightarrow H, \qquad (t,x) \longmapsto b(t,x) \]
is continuous. Note that

\customlabel{K-C}
\[ (\mu_0+\mu_h)[K^c] = \underbrace{\mu_0[K^c]}_{\leq\bar\delta\leq\frac\eps2} + \mu_h[K^c] \leq \frac\eps2 + \!\!\!\! \underbrace{\int\limits_{K^c} \rho \d \mu_0(t, x)}_{\leq\frac\eps2 \text{ by } \eqref{TH2-BARDELTA} \text{ and } \eqref{TH2-RHO}} \!\!\!\! \leq \eps .  \tag{\ref{K-C}} \]

Applying Dugundji's Extension Theorem (see \cite[Theorem 4.1]{Dug51}) to the function $b|_K$ guarantees that there exists a continuous function $\bar b \colon [0,1] \times H \longrightarrow H$ with $\|\bar b\|_\infty \leq \|b\|_\infty$ and $\|\bar b\|_{\infty,A} \leq \|b\|_{\infty,A}$ which coincides with $b$ on $K$. Starting from the left-hand side of the assertion we have

\[ \E \exp \frac{\beta_A}{\|h\|_{\infty}^2} \left| \int\limits_0^1 b(t, Z_t^{\ell A} + h(t)) - b(t, Z_t^{\ell A}) \d t \right|_{H}^2  .  \]

Adding and subtracting $\bar b$ and using that $b-\bar b=0$ on $K$ yields that the above equals

\[ \E \exp \frac{\beta_A}{\|h\|_{\infty}^2} \left| \int\limits_0^1 \bbm 1_{K^c}(t, Z_t^{\ell A} + h(t))  \left[ \smash{\underbrace{b(t, Z_t^{\ell A} + h(t)) - \bar b(t, Z_t^{\ell A} + h(t))}_{\in [-2,2]}} \right] \right. \]
\[  \hspace{-6mm} \left. - \bbm 1_{K^c}(t, Z_t^{\ell A}) \left[ \smash{\underbrace{b(t, Z_t^{\ell A}) - \bar b(t, Z_t^{\ell A})}_{\in [-2,2]}} \right]  \d t  \right.  \]
\[ \hspace{-6mm} \left.   +   \int\limits_0^1 \bar b(t, Z_t^{\ell A} + h(t)) - \bar b(t, Z_t^{\ell A}) \d t \right|_{H}^2 . \]

Applying the fact that $(a+b)^2 \leq 2a^2 + 2b^2$ we estimate from above by

\[   \E \exp \left( \frac{8 \beta_A}{\|h\|_{\infty}^2} \left( \int\limits_0^1 \bbm 1_{K^c}(t, Z_t^{\ell A} + h(t))  +  \bbm 1_{K^c}(t, Z_t^{\ell A}) \d t \right)^2 \right. \]
\[   \hspace{1mm}  +   \frac{2 \beta_A}{\|h\|_{\infty}^2} \left.\left| \int\limits_0^1  \bar b(t, Z_t^{\ell A} + h(t)) - \bar b(t, Z_t^{\ell A}) \d t \right|_{H}^2 \right) \]

\[ \hspace{-1mm} = \E \exp \left( \frac{8 \beta_A}{\|h\|_{\infty}^2} \left( \int\limits_0^1 \bbm 1_{K^c}(t, Z_t^{\ell A} + h(t)) + \bbm 1_{K^c}(t, Z_t^{\ell A} ) \d t \right)^2  \right) \]
\[  \hspace{-8mm}  \cdot  \exp \left( \frac{2 \beta_A}{\|h\|_{\infty}^2} \left| \int\limits_0^1 \bar b(t, Z_t^{\ell A} + h(t)) - \bar b(t, Z_t^{\ell A}) \d t \right|_{H}^2 \right) \]

and using Young's Inequality this is bounded by
\[ \frac12 \underbrace{\E \exp \left( \frac{16 \beta_A}{\|h\|_{\infty}^2} \left| \int\limits_0^1 \bbm 1_{K^c}(t, Z_t^{\ell A} + h(t)) + \bbm 1_{K^c}(t, Z_t^{\ell A})  \d t \right|_H^2  \right)}_{=:\ \!A_1}   \]
\[  \hspace{-8mm}  +  \frac12 \underbrace{\E \exp \left( \frac{4 \beta_A}{\|h\|_{\infty}^2} \left| \int\limits_0^1 \bar b(t, Z_t^{\ell A} + h(t)) - \bar b(t, Z_t^{\ell A}) \d t \right|_{H}^2 \right)}_{=:\ \!A_2} . \]

Let us estimate $A_1$ first

\[ A_1 = 1 + \sum\limits_{k=1}^\infty \frac1{k!} \left(\frac{16 \beta_A}{\|h\|_{\infty}^2}\right)^k \E \left| \int\limits_0^1 \bbm 1_{K^c}(t, Z_t^{\ell A} + h(t)) + \bbm 1_{K^c}(t, Z_t^{\ell A})  \d t \right|_H^{2k} \]

\[ \leq 1 + \sum\limits_{k=1}^\infty \frac1{k!} \left(\frac{16 \beta_A}{\|h\|_{\infty}^2}\right)^k 2^{2k} \underbrace{\left(\mu_h[K^c] + \mu_0[K^c]\right)}_{\leq \eps \text{ by } \eqref{K-C}} \leq 1 + \sum\limits_{k=1}^\infty \frac1{k!} \left(\frac{64 \beta_A}{\|h\|_{\infty}^2}\right)^k \eps \]
\[ \leq 1 + \exp\left( \frac{64 \beta_A}{\|h\|_{\infty}^2} \right) \eps = 1 + 1 = 2 . \]

This concludes the estimate for $A_1$. Let us now estimate $A_2$. Since $\bar b$ is continuous there exists a sequence $\bar b^{(m)} \colon [0,1] \times H \longrightarrow H$ of functions with $\|\bar b^{(m)}\|_\infty \leq 1$ and $\|\bar b^{(m)}\|_{\infty,A} \leq 1$ which are smooth in the second component (i.e.~twice continuously differentiable) such that $\bar b^{(m)}$ converges to $\bar b$ everywhere, i.e.
\[ \hspace{10mm}  \bar b^{(m)}(t,x) \overset{m\rightarrow\infty}\longrightarrow \bar b(t,x),   \qquad\quad \forall t\in[0,1] , \ \forall x\in H . \]
Using the above considerations $A_2$ equals
\[  \E  \exp \frac{4 \beta_A}{\|h\|_{\infty}^2}   \left| \int\limits_0^1   \lim\limits_{m\rightarrow\infty}   \bar b^{(m)} (t, Z^{\ell A}_t + h(t))  -  \bar b^{(m)} (t, Z^{\ell A}_t ) \d t \right|_H^2 , \]

which in turn can be bounded using Fatou's Lemma by
\customlabel{LIMINF}
\[  \liminf\limits_{m\rightarrow\infty} \E  \exp \frac{4 \beta_A}{\|h\|_{\infty}^2}   \left| \int\limits_0^1     \bar b^{(m)} (t, Z_t^{\ell A} + h(t))  -  \bar b^{(m)} (t, Z_t^{\ell A} ) \d t \right|_H^2 .  \tag{\ref{LIMINF}}  \]

Applying Step $1$ with $b$ replaced by $\bar b^{(m)}$ yields that \eqref{LIMINF} and henceforth $A_2$ is bounded by $C$, so that in conclusion we have

\[ \E \exp \frac{\beta_A}{\|h\|_{\infty}^2} \left| \int\limits_0^1 b(t, Z_t^{\ell A} + h(t)) - b(t, Z_t^{\ell A}) \d t \right|_{H}^2  \leq \frac12 A_1 + \frac12 A_2 \leq 1 + \frac C2 \leq 3 , \]
which completes the proof.

\qed
\end{prof}

\section{A concentration of measure result}

For this section let us define

\[ \hspace{9mm}  Z^A(t,x) := Z^A_t + e^{-tA} x,  \hspace{4mm} \qquad \forall x \in H, \ t\in[0,\infty[ \]

then for every $x\in H$, $Z(\cdot,x)$ is an Ornstein--Uhlenbeck process starting in $x$. Furthermore, we define the image measure
\[ \P_x := \P \circ Z^A(\cdot, x)^{-1}, \qquad \forall x\in H \]
and the projections
\[ \pi_t(f) := f(t), \qquad \forall f\in \ca C( [0,\infty[, H), \ t\in[0,\infty[ ,  \]

which come with their canonical filtration

\[ \bar{\ca G}_t := \sigma ( \pi_s | s \leq t ) , \]

so that $(\ca C( [0,\infty[, H), (\P_x)_{x\in H}, (\pi_t)_{t\in[0,\infty[}, (\bar{\ca G}_t)_{t\in[0,\infty[})$ is a universal Markov process (see \cite[Proposition 4.3.5]{LR15} and \cite[Section 42]{Bau96} or \cite[Section 3.4]{Jac05} for the definition of a universal Markov process). Additionally, we set

\[ \ca G_t := \{ Z^{-1}(B) | B \in \bar{\ca G}_t \} \]

as the initial sigma-algebra, so that $Z$ becomes $\ca G_t / \bar{\ca G}_t$-measurable.

\begin{cor}
\label{THM2}
\customtoc{Corollary}

There exists $\beta_A>0$ (depending only on the drift term $A$ of the Ornstein--Uhlenbeck process $Z^A$) such that for all $0\leq r<u\leq1$ and for any bounded Borel measurable function $b\colon [r,u]\times H \longrightarrow H$ with $\|b\|_\infty \leq 1$, $\|b\|_{\infty,A} \leq 1$, any bounded Borel measurable functions $h_1, h_2\colon [r,u]  \longrightarrow H$ with
\[ \sum\limits_{n\in\N} |h_{1,n}(t)|^2 \lambda_n^2  +  \sum\limits_{n\in\N} |h_{2,n}(t)|^2 \lambda_n^2 < \infty ,  \qquad \forall t \in [0,1]  . \]
for any $\eta \geq 0$ the inequality
\[ \P \left[ \left. \left| \int\limits_r^u b(s, Z_s^A + h_1(s)) - b(s, Z_s^A + h_2(s)) \d s \right|_H  >  \eta \ell^{1/2} \|h_1-h_2\|_\infty \right| \ca G_r \right] \leq 3 e^{- \beta_A \eta^2} \]
holds, where $\ell:=u-r$.
\end{cor}

\begin{prof}

Let $r,u,\ell,b,h_1$ and $h_2$ be as in the assertion. Note that the assertion is trivial if $h_1=h_2$, hence w.lo.g.~we assume $\|h_1 - h_2\|_\infty \neq 0$. We define the stochastic processes $\tilde Z^{\ell A}_t := \ell^{-1/2} Z_{\ell t}^A$ and $\tilde B_t := \ell^{-1/2} B_{\ell t}$. Note that $\tilde B$ is again a Brownian motion w.r.t.~the normal, right-continuous filtration $(\tilde{\ca F}^\ell_t)_{t\in[0,\infty[} := (\ca F_{\ell t})_{t\in[0,\infty[}$. Additionally, we have
\[ \hspace{-25mm}  \tilde Z^{\ell A}_t = \ell^{-1/2} Z_{\ell t}^A = \ell^{-1/2} \int\limits_0^{\ell t} e^{(\ell t - s) A} \d B_s \]
\[ = \int\limits_0^{\ell t} e^{\ell \left(t - \frac s\ell\right) A} \ell^{-1/2} \d B_{\frac s\ell \ell} = \int\limits_0^t e^{(t-s') \ell A} \d \tilde B_{s'} . \]
Hence, $\tilde Z^{\ell A}$ is an Ornstein--Uhlenbeck process with drift term $\ell A$.

For the reader's convenience we add the integration variable as a superscript to the respective measure which we integrate against, hence the left-hand side of the claim reads
\[ \P^{(\mrm d\omega)} \left[ \left. \left| \int\limits_r^u b(s, Z_s^A(\omega) + h_1(s)) - b(s, Z_s^A(\omega) + h_2(s)) \d s \right|_H  >  \eta \ell^{1/2} \|h_1-h_2\|_\infty \right| \ca G_r \right] . \]

Fix an $\omega'\in\Omega$. Using the transformation $s' := \ell^{-1}(s-r)$ this equals

\[ \hspace{-52mm}  \P^{(\mrm d\omega)} \left[ \left. \ell \left| \int\limits_0^1 b(\ell s'+r, Z_{\ell s'+r}^A(\omega) + h_1(\ell s'+r)) \right. \right. \right. \]
\[ \hspace{25mm}  \left.\left.\left. \vphantom{\int\limits_0^1} - b(\ell s'+r, Z_{\ell s'+r}^A(\omega) + h_2(\ell s'+r)) \d s' \right|_{H}  >  \eta \ell^{1/2} \|h_1-h_2\|_\infty   \right| \ca G_r  \right] (\omega') .  \]

Recall the definitions of $\pi_t$ and $\bar{\ca G_t}$ at the beginning of this section. Since $\ca G_t$ is the initial sigma-algebra of $\bar{\ca G_t}$ w.r.t.~$Z^A$ we have

\[ \E \left[ \pi_t \circ Z^A | \ca G_r \right] (\omega')  =  \E_0 \left[ \pi_t  | \bar{\ca G_r} \right] (Z^A(\omega')) ,  \]

where $\E_0$ denotes the expectation w.r.t.~the measure $\P_0$. Applying this to the above situation we obtain that the left-hand side of the assertion reads

\[ \hspace{-64mm}  \P_0^{(\mrm d\omega)} \left[ \left| \int\limits_0^1 b(\ell s+r, \pi_{\ell s+r}(\omega) + h_1(\ell s+r)) \right. \right. \]
\[ \hspace{20mm}  \left.\left. \left. \vphantom{\int\limits_0^1} - b(\ell s+r, \pi_{\ell s+r}(\omega) + h_2(\ell s+r)) \d s \right|_H  >  \eta \ell^{-1/2} \|h_1-h_2\|_\infty    \right| \bar{\ca G}_r \right] (Z^A(\omega')) , \]

Applying the universal Markov property (see \cite[Equation (42.18)]{Bau96} or \cite[Equation (3.108)]{Jac05}) we have

\customlabel{TH2EQ1}
\[ \hspace{-64mm} = \P_{\pi_r(Z^A(\omega'))}^{(\mrm d\omega)} \left[ \left| \int\limits_0^1 b(\ell s+r, \pi_{\ell s}(\omega) + h_1(\ell s+r)) \right. \right. \]
\[ \hspace{32mm} \left. \left. \vphantom{\int\limits_0^1}  - b(\ell s+r, \pi_{\ell s}(\omega) + h_2(\ell s+r)) \d s \right|_H  > \eta \ell^{-1/2} \|h_1-h_2\|_{\infty} \right]  . \tag{\ref{TH2EQ1}} \]

We define
\[ \tilde b(t,x)   :=             b  (\ell t+r, \ell^{1/2} x), \]
\[ \tilde h_1(t)   := \ell^{-1/2} h_1(\ell t+r), \]
\[ \tilde h_2(t)   := \ell^{-1/2} h_2(\ell t+r), \]

so that the Expression \eqref{TH2EQ1} simplifies to
\[ \P_{\pi_r(Z^A(\omega'))}^{(\mrm d\omega)} \left[ \left| \int\limits_0^1 \tilde b(s, \ell^{-1/2} \pi_{\ell s} + \tilde h_1(s)) - \tilde b(s, \ell^{-1/2} \pi_{\ell s} + \tilde h_2(s)) \d s \right|_H > \eta \left\|\tilde h_1 - \tilde h_2 \right\|_\infty \right] .  \]

Note that $\tilde b$, $\tilde h_1$, $\tilde h_2$ are all bounded Borel measurable functions and $\| \tilde b \|_\infty = \| b \|_\infty \leq 1$ as well as $\| \tilde b \|_{\infty,A} = \| b \|_{\infty,A} \leq 1$ hold. Plugging in the definition of $\P_x$ the above reads

\[ \hspace{-54mm}  (\P \circ Z^A(\cdot, Z^A_r(\omega'))^{-1} )^{(\mrm d\omega)} \left[ \left| \int\limits_0^1 \tilde b (s, \ell^{-1/2} \pi_{\ell s}(\omega) + \tilde h_1(s)) \right.\right. \]
\[ \hspace{35mm}  \left.\left. \vphantom{\int\limits_0^1} - \tilde b(s, \ell^{-1/2} \pi_{\ell s}(\omega) + \tilde h_2(s)) \d s \right|_{H}  >  \eta \left\| \tilde h_1 - \tilde h_2 \right\|_\infty \right]   \]

\[ \hspace{-17mm}  = \P \left[ \left| \int\limits_0^1 \tilde b_{\omega',\tilde h_2} (s, \underbrace{\ell^{-1/2} Z^A(\ell s, Z^A_r(\omega')) - \ell^{-1/2} e^{-\ell sA} Z^A_r(\omega')}_{= \ell^{-1/2} Z^A_{\ell s} = \tilde Z^{\ell A}_s} + \tilde h_1(s) - \tilde h_2(s)) \right. \right. \]
\[ \hspace{9mm}  \left. \left. \vphantom{\int\limits_0^1} - \tilde b_{\omega',\tilde h_2}(s, \underbrace{\ell^{-1/2} Z^A(\ell s, Z^A_r(\omega')) - \ell^{-1/2} e^{-\ell sA} Z^A_r(\omega')}_{= \ell^{-1/2} Z^A_{\ell s} = \tilde Z^{\ell A}_s}) \d s \right|_H  >  \eta \left\| \tilde h_1 - \tilde h_2 \right\|_\infty \right]  , \]

where $\tilde b_{\omega',\tilde h_2}(t,x) := \tilde b(t,x + \ell^{-1/2} e^{-\ell tA} Z^A_r(\omega') + \tilde h_2(t))$. Recall that $\tilde Z^{\ell A}$ is an Ornstein--Uhlenbeck process which starts in $0$. By Theorem \ref{THM1B} there exist constants $\beta_A$ (depending on the drift term $A$, but independent of $\ell$ since $\ell\in \ ]0,1]$) and $C$ such that the conclusion of Theorem \ref{THM1B} holds for every Ornstein--Uhlenbeck process $\tilde Z^{\ell A}$ with the same constants $\beta_A$ and $C$. Since $\exp(\beta_A |\cdot|^2)$ is increasing on $\R_+$ the above equals

\[ \P \left[ \exp \left( \frac{\beta_A}{\left\|\tilde h_1 - \tilde h_2 \right\|^2_{\infty}} \left| \int\limits_0^1 \tilde b_{\omega',\tilde h_2}(s, \tilde Z^{\ell A}_s + \tilde h_1(s) - \tilde h_2(s)) - \tilde b_{\omega',\tilde h_2}(s, \tilde Z^{\ell A}_s ) \d s \right|^2_H \right)  >  \exp\left( \beta_A \eta^2 \right) \right] \]

and by Chebyshev's Inequality this can be estimated from above by
\[ e^{- \beta_A \eta^2} \E \exp \left( \frac{\beta_A}{\left\| \tilde h_1 - \tilde h_2 \right\|_{\infty}^2} \left| \int\limits_0^1  \tilde b_{\omega',\tilde h_2}(s, \tilde Z^{\ell A}_s + \tilde h_1(s) - \tilde h_2(s)) - \tilde b_{\omega',\tilde h_2}(s, \tilde Z^{\ell A}_s) \d s \right|^2_H \right) . \]

Since $\| \tilde b_{\omega',\tilde h_2} \|_\infty = \| \tilde b \|_\infty \leq 1$ as well as $\| \tilde b_{\omega',\tilde h_2} \|_{\infty,A} = \| \tilde b \|_{\infty,A} \leq 1$ holds the conclusion of Theorem \ref{THM1B} implies that the above expression is smaller than
\[ C e^{- \beta_A \eta^2} \leq 3 e^{- \beta_A \eta^2} .  \]
\qed
\end{prof}

\begin{cor}
\label{THM4}
\customtoc{Corollary}

For all $0\leq r<u\leq1$ and for every bounded Borel measurable function $b\colon [r,u]\times H \longrightarrow H$ with $\|b\|_\infty \leq 1$, $\|b\|_{\infty,A} \leq 1$ and for all bounded $\ca G_r$-measurable random variables $x,y \colon \Omega \longrightarrow H$ with
\[ \sum\limits_{n\in\N} |x_n|^2 \lambda_n^2  +  \sum\limits_{n\in\N} |y_n|^2 \lambda_n^2 < \infty  . \]
We have
\[ \E \left[ \left| \int\limits_r^u b(s, Z_s^A + x) - b(s, Z_s^A + y) \d s \right|_{H}^p \left. \vphantom{\int\limits_r^u} \right| \ca G_r \right]  \leq  3 \beta_A^{p/2} p^{p/2} \ell^{p/2} |x-y|_H^p  , \]
where $\ell:=u-r$, $p\in\N$ and $\beta_A > 0$ is the constant from Corollary \ref{THM2}.
\end{cor}

\begin{prof}

Let $0\leq r<u\leq1$ and $b,p,\ell$ as in the assertion.\\

\textbf{Step 1:} Deterministic $x$, $y$\\

Let $x,y \in H$ be non-random with $x \neq y$. We set

\[ S := \beta_A^{1/2} \ell^{-1/2} |x-y|_H^{-1} \left| \int\limits_r^u b(s, Z_s^A + x) - b(s, Z_s^A + y) \d s \right|_H  \]

and calculate

\[ \E \left[ \left. S^p \right| \ca G_r \right]   =   \E \left[ \left.   \int\limits_0^\infty  \bbm 1_{\{S > \eta\}} p \eta^{p-1} \d \eta \right| \ca G_r \right] .  \]

Notice that the above is valid since $S$ is a non-negative random variable. Using Fubini's Theorem the above equals

\[ \int\limits_0^\infty p \eta^{p-1} \P \left[ \left.  S > \eta  \right| \ca G_r \right] \d \eta .  \]

Plugging in the definition of $S$ the above line reads

\[ \int\limits_0^\infty p \eta^{p-1} \P \left[ \left.  \left| \int\limits_r^u b(s, Z_s^A + x) - b(s, Z_s^A + y) \d s \right|_H  \! > \beta_A^{-1/2} \eta \ell^{1/2} |x-y|_H \right| \ca G_r \right] \d \eta .  \]

We estimate the probability inside the integral by invoking Corollary \ref{THM2}, hence the above expression is smaller than

\[ 3 \int\limits_0^\infty p \eta^{p-1} e^{-\eta^2} \d \eta  =  \frac{3p}2 \Gamma\left( \frac p2 \right) .  \]

Using Stirling's formula this is bounded from above by

\[ \frac32 \underbrace{p \sqrt{\frac{4\pi}p} 2^{-p/2} e^{-p/2} e^{\frac1{6p}}}_{\leq \sqrt{2\pi} e^{-1/2} e^{\frac16}} p^{p/2} \leq 3 p^{p/2} ,  \]

which proves that $\E [ S^p | \ca G_r ] \leq 3 \beta_A^{p/2} p^{p/2}$, concluding the assertion in the case that $x$ and $y$ are deterministic.\\

\textbf{Step 2:} Random $x$, $y$\\

Let $x,y\colon \Omega \longrightarrow H$ be $\ca G_r$ measurable random variables of the form

\[ x = \sum\limits_{i=1}^n \bbm 1_{A_i} x_i ,  \qquad y = \sum\limits_{i=1}^n \bbm 1_{A_i} y_i ,  \]

where $x_i,y_i\in H$ and $(A_i)_{1\leq i \leq n}$ are pairwise disjoint sets in $\ca G_r$. Notice that due to the disjointness we have

\[ b\left(t, Z_t^A  + \sum\limits_{i=1}^n \bbm 1_{A_i} x_i \right) - b\left(t, Z_t^A + \sum\limits_{i=1}^n \bbm 1_{A_i} y_i \right)  =  \sum\limits_{i=1}^n \bbm 1_{A_i} \left[ b(t, Z_t^A  +  x_i ) - b(t, Z_t^A + y_i) \right] .   \]

Let $p$ be a positive integer. Starting from the left-hand side of the assertion an using the above identity yields

\[  \E \left[ \left. \left|  \int\limits_r^u b(t, Z_t^A  + x) - b(t, Z_t^A + y) \d t \right|_H^p  \right| \ca G_r \right]  \]

\[ =  \sum\limits_{i=1}^n \E \left[ \left. \bbm 1_{A_i} \left| \int\limits_r^u b(t, Z_t^A  +  x_i) - b(t, Z_t^A + y_i) \d t \right|_H^p  \right| \ca G_r \right] .  \]

Since $A_i \in \ca G_r$ this can be expressed as

\[ \sum\limits_{i=1}^n \bbm 1_{A_i} \E\left[ \left. \left| \int\limits_r^u b(t, Z_t^A + x_i) - b(t, Z_t^A + y_i) \d t \right|_H^p \right| \ca G_r \right]  \]

and by invoking Step 1 this is bounded from above by

\[ 3 \beta_A^{p/2} p^{p/2} \ell^{p/2} \sum\limits_{i=1}^n \bbm 1_{A_i} |x_i-y_i|_H^p   =  3 \beta_A^{p/2} p^{p/2} \ell^{p/2} | x-y |_H^p .  \]

In conclusion we obtained the result for step functions $x$, $y$.
The result for general $\ca G_r$ measurable random variables $x,y$ now follows by approximation via step functions and taking limits.

\qed

\end{prof}


\begin{footnotesize}

\end{footnotesize}


\begin{thebibliography}{MNS89}

\bibitem[Bau96]{Bau96}
H.~Bauer.
\newblock {\em Probability Theory}.
\newblock De Gruyter Studies in Mathematics \textbf{23}, pp.~523, 1996.

\bibitem[BJ97]{BJ97}
X.~Bardina, M.~Jolis.
\newblock {\em{An extension of Ito's formula for elliptic diffusion
  processes}}.
\newblock { {Stoch.~Proc.~Appl.}} \textbf{69}, no.~1, pp.~83--109, 1997.

\bibitem[Bog98]{Bog98}
V.~I.~Boga{\v{c}}ev.
\newblock {\em Gaussian measures}, { Mathematical surveys and
 monographs} \textbf{62}.
\newblock Am.~Math.~Soc., 433 pages, 1998.

\bibitem[BR07]{BR07}
X.~Bardina, C.~Rovira.
\newblock {\em On It{\^o}'s formula for elliptic diffusion processes}.
\newblock {Bernoulli}, \textbf{13} no.~3, pp.~820--830, 2007.
\newblock \href{http://arxiv.org/abs/0709.0627}{arXiv: 0709.0627v1}.

\bibitem[Dav76]{Dav76}
B.~Davis.
\newblock {\em On the $l^p$ norms of stochastic integrals and other martingales}.
\newblock {Duke Math.~J.}, \textbf{43} no.~4, pp.~697--704, 1976.

\bibitem[Dug51]{Dug51}
J.~Dugundji.
\newblock {\em An extension of Tietze's theorem}.
\newblock {Pac.~J.~Math.}, \textbf{1}, no.~3, pp.~353--367, 1951.

\bibitem[DZ92]{DZ92}
G.~Da Prato, J.~Zabczyk.
\newblock{\em Stochastic Equations in Infinite Dimensions}.
\newblock{Cambridge Univ. Pr.},
pp.~454, 1992.

\bibitem[Jac05]{Jac05}
N.~Jacob.
\newblock {\em Pseudo Differential Operators \& Markov Processes: Markov processes and applications}.
\newblock Imperial College Press, pp.~440, 2005.

\bibitem[Kho14]{Kho14}
D.~Khoshnevisan.
\newblock{\em Analysis of stochastic partial differential equations}
\newblock{CBMS regional conference series in mathematics} \textbf{119}.
\newblock Am.~Math.~Soc., 2014.

\bibitem[LR15]{LR15}
W.~Liu, M.~R{\"o}ckner.
\newblock {\em Stochastic Partial Differential Equations: An Introduction}.
\newblock Springer International Publishing, pp.~266, 2015.

\bibitem[MNS89]{MNS89}
A.~Millet, D.~Nualart, M.~Sanz.
\newblock{\em{Integration by Parts and Time Reversal for Diffusion
  Processes}}.
\newblock{Ann.~Probab.}, \textbf{17} no.~1, pp.~208--238, 1989.

\bibitem[{\O}ks10]{Oks10}
B.~{\O}ksendal.
\newblock{\em {Stochastic Differential Equations: An Introduction with
  Applications}}.
\newblock Springer, 6th ed.~2003, corr.~6th printing 2014 edition, 2010.

\bibitem[Ose12]{Ose12}
A.~Osekowski.
\newblock {\em Sharp Martingale and Semimartingale Inequalities},
{ Monografie Matematyczne} \textbf{72}
\newblock Springer, 2012.

\bibitem[Par86]{Par86}
E.~Pardoux.
\newblock {\em Grossissement d'une filtration et retournement du temps d'une
  diffusion. (french) [enlargement of a filtration and time reversal of a
  diffusion]}.
\newblock { S\'eminaire de Probabilit\'es XX, 1984/85}, Lecture notes in
  Mathematics \textbf{1204}, pp.~ 48--55. Springer, 1986.

\bibitem[RSZ08]{RSZ08}
M.~R{\"o}ckner, B.~Schmuland, X.~Zhang.
\newblock {\em Yamada--{W}atanabe {T}heorem for stochastic evolution equations in infinite dimensions}.
\newblock { Cond.~Matt.~Phys.}, \textbf{11} no.~2, pp.~247--259, 2008.

\bibitem[{Sha}14]{Sha14}
A.~{Shaposhnikov}.
\newblock {\em{Some remarks on Davie's uniqueness theorem}}.
\newblock { ArXiv e-prints}, 2014.
\newblock \href{http://arxiv.org/abs/1401.5455}{arXiv: 1401.5455}.

\bibitem[Tao11]{Tao11}
T.~Tao.
\newblock {\em An introduction to measure theory}.
\newblock Graduate studies in mathematics \textbf{126}. Am.~Math.~Soc., 2011.

\bibitem[Wan91]{Wan91}
G.~Wang.
\newblock {\em Some sharp inequalities for conditionally symmetric martingales}.
\newblock {T.~Am.~Math.~Soc.},
  \textbf{328}, pp.~392--419, 1991.

\end{thebibliography}
\end{document}